# A stochastic fixed point equation for weighted minima and maxima

Gerold Alsmeyer[a] and Uwe Rösler[b]

[a]Institut für Mathematische Statistik, Fachbereich Mathematik, Westfälische Wilhelms-Universität Münster, Einsteinstraße 62, D-48149 Münster, Germany
[b]Mathematisches Seminar, Christian-Albrechts-Universität Kiel, Ludewig-Meyn-Straße 4, D-24098 Kiel, Germany

**Abstract.** Given any finite or countable collection of real numbers $T_j, j \in J$, we find all solutions $F$ to the stochastic fixed point equation

$$W \stackrel{\mathrm{d}}{=} \inf_{j \in J} T_j W_j,$$

where $W$ and the $W_j, j \in J$, are independent real-valued random variables with distribution $F$ and $\stackrel{\mathrm{d}}{=}$ means equality in distribution. The bulk of the necessary analysis is spent on the case when $|J| \geq 2$ and all $T_j$ are (strictly) positive. Nontrivial solutions are then concentrated on either the positive or negative half line. In the most interesting (and difficult) situation $T$ has a characteristic exponent $\alpha$ given by $\sum_{j \in J} T_j^\alpha = 1$ and the set of solutions depends on the closed multiplicative subgroup of $\mathbb{R}^{>} = (0, \infty)$ generated by the $T_j$ which is either $\{1\}$, $\mathbb{R}^{>}$ itself or $r^{\mathbb{Z}} = \{r^n \colon n \in \mathbb{Z}\}$ for some $r > 1$. The first case being trivial, the nontrivial fixed points in the second case are either Weibull distributions or their reciprocal reflections to the negative half line (when represented by random variables), while in the third case further periodic solutions arise. Our analysis builds on the observation that the logarithmic survival function of any fixed point is harmonic with respect to $\Lambda = \sum_{j \geq 1} \delta_{T_j}$, i.e. $\Gamma = \Gamma \star \Lambda$, where $\star$ means multiplicative convolution. This will enable us to apply the powerful Choquet–Deny theorem.

**Résumé.** Étant donné un ensemble fini ou dénombrable de nombres réel $T_j, j \in J$, nous trouvons l'ensemble des solutions $F$ de l'équation fonctionelle

$$W \stackrel{\mathrm{d}}{=} \inf_{j \in J} T_j W_j,$$

où $W$ et les $W_j, j \in J$, sont des variables aléatoires mutuellement indépendantes ayant la loi $F$ et $\stackrel{\mathrm{d}}{=}$ signifie identité en loi. L'essentiel de ce travail concerne le cas où $|J| \geq 2$ et tous les $T_j$ sont (strictement) positifs. Dans ce cas, toutes les solutions sont concentrées soit sur $(0, \infty)$ soit sur $(-\infty, 0)$. Dans la situation la plus intéressante (et plus difficile) $T$ a un exposant charactéristique $\alpha$ donné par $\sum_{j \in J} T_j^\alpha = 1$, et l'ensemble des solutions dépend du sous-groupe multiplicatif de $\mathbb{R}^{>} = (0, \infty)$ généré par les $T_j$, qui est $\{1\}, \mathbb{R}^{>}$ lui-meme, ou $r^{\mathbb{Z}} = \{r^n \colon n \in \mathbb{Z}\}$ pour quelque $r > 1$. Le premier cas etant trivial, les points fixes non-triviaux dans le second cas sont ou bien les lois de Weibull ou bien leurs images réciproques sur $(-\infty, 0)$ (si elles sont représentées par des variables aléatoires). Dans le troisième cas, il y a des solutions périodiques supplémentaires. Notre analyse est basée sur l'observation que le logarithme de la fonction de survie de chaque point fixe est harmonique relatif à $\Lambda = \sum_{j \geq 1} \delta_{T_j}$, c'est-à-dire $\Gamma = \Gamma \star \Lambda$, où $\star$ dénote la convolution multiplicative. Cela nous permettrons l'utilisation du theorème puissant de Choquet et Deny.

## 1. Introduction

Given a finite or infinite sequence $T = (T_j)_{j \in J}$ of real numbers (thus $J = \{1, \ldots, n\}$ or $J = \mathbb{N}$), consider the stochastic fixed point equations

$$W \stackrel{d}{=} \inf_{j \in J} T_j W_j \tag{1.1}$$

and

$$W \stackrel{d}{=} \sup_{j \in J} T_j W_j, \tag{1.2}$$

for i.i.d. real-valued random variables $W, W_1, W_2, \ldots$. The general goal is to determine the collections $\mathfrak{F}_T^{\min}$ and $\mathfrak{F}_T^{\max}$ of all (distributional) fixed points of these equations, that is all distributions of $W$ such that (1.1), respectively (1.2) holds true. For the more general situation of random weights $T_j$, these equations are discussed in some detail by Jagers and the second author [6], while Neininger and Rüschendorf [7] provide examples from the asymptotic analysis of recursive algorithms and data structures where solutions of equations of this type emerge as limiting distributions. More examples from various areas of applied probability which are related to max- or min-type fixed point equations appear in a recent long survey paper by Aldous and Bandyopadhyay [1]. These include, for instance, the extinction time of Galton–Watson processes (= height of Galton–Watson trees) and the extremal positions in as well as the range of branching random walks.

To our best knowledge the problem of providing all fixed points of (1.1) or (1.2) for general random $T_j$ is a completely open one, and the present article contributes to it by giving a complete answer for the simpler case of deterministic $T_j$. Our main motivation for a detailed analysis of this special case, though of interest in its own right, is to learn about how to tackle the general situation. Work in progress gives rise to the conjecture that the solutions of (1.1) or (1.2) for random $T_j$ may be described as suitable mixtures of those for deterministic $T_j$ where the mixing distribution is itself a solution to a related fixed point equation. We refer to a future publication.

The trivial equivalence

$$\mathcal{L}(W) \in \mathfrak{F}_T^{\max} \quad \Longleftrightarrow \quad \mathcal{L}(-W) \in \mathfrak{F}_T^{\min} \tag{1.3}$$

shows that we must only analyze (1.1). Here $\mathcal{L}(X)$ denotes the distribution of a random variable $X$. It further suffices to consider the case where all weights $T_j, j \in J$, are nonzero because then the fixed points of the equation

$$W \stackrel{d}{=} 0 \wedge \inf_{j \in J} T_j W_j, \tag{1.4}$$

that is of (1.1) extended by a zero weight, are just those of (1.1) concentrated on $(-\infty, 0]$. Direct inspection shows that only $F = \delta_0$, the Dirac measure at 0, solves (1.4) if all $T_j$ are negative. Hence the nontrivial analysis of (1.1) reduces to the following three cases:

(C1) all $T_j$, $j \in J$, are positive ($> 0$);
(C2) all $T_j$, $j \in J$, are negative;

(C3) both, $J^> \stackrel{\text{def}}{=} \{j \in J: T_j > 0\}$ and $J^< \stackrel{\text{def}}{=} \{j \in J: T_j < 0\}$ are nonempty.

Cases (C2) and (C3) will be dealt with rather shortly in the final section, the latter by drawing on the results for (C1) in an appropriate manner. As for (C2), we will obtain that $\mathfrak{F}_T = \{\delta_0\}$ if $J = \mathbb{N}$, while for finite $J$ our result only provides a characterization of the solutions which relates them to another stochastic fixed point equation of minimax-type (see (5.1) and (5.3)). In a special case, namely $T_1 = \cdots = T_n = c < -1$, this minimax-type equation also shows up in the context of randomized game tree evaluation for which Khan, Devroye and Neininger [2] could show that a nontrivial fixed point exists. Recent work of the first author with Meiners [3] on the very same equation but for general $T_j$ even shows that the solution set can be quite large. However, his methods are very different from those employed here. For now we are thus left with the case (C1) and therefore assume until further notice that *all $T_j$ are positive.*

In the following we will always write $\mathfrak{F}_T$ instead of $\mathfrak{F}_T^{\min}$. For our convenience, it is also stipulated as quite common in probability theory that the same symbol $F$ is used for a distribution on $\mathbb{R}$ as well as for its left continuous distribution function, so $F(t)$ means the same as $F((-\infty, t))$. In order to gain further insight into the posed problem and to provide an outline of the necessary analysis we begin with some basic observations the simplest one being that $\mathfrak{F}_T$ always contains the trivial solution $\delta_0$. If $|J| = 1$ there is clearly no other fixed point unless $T_1 = 1$. In the latter case Eq. (1.1) becomes trivial and $\mathfrak{F}_T$ consists in fact of all distributions on $\mathbb{R}$. We hence make the *standing assumption*:

$$|J| \geq 2$$

hereafter. A fixed point is called nontrivial if it is not a Dirac measure $\delta_c$ for some $c \in \mathbb{R}$.

A good starting point made upon direct inspection is that the exponential distributions $\text{Exp}(c), c > 0$, are nontrivial fixed points whenever $\sum_{j \in J} T_j^{-1} = 1$. We will show in Section 4 of this article that they are then indeed the only ones if the multiplicative group generated by the $T_j$ (denoted as $\mathbb{G}_\Lambda$ below) is $\mathbb{R}^> \stackrel{\text{def}}{=} (0, \infty)$, whereas further periodic solutions exist otherwise. They will also be defined there. With this result at hand the situation where

$$\sum_{j \in J} T_j^{-\alpha} = 1 \tag{1.5}$$

for some (necessarily unique) $\alpha > 0$ is also settled because (1.1) is equivalent to

$$W^\alpha \stackrel{\text{d}}{=} \inf_{j \in J} T_j^\alpha W_j^\alpha \tag{1.6}$$

for any $\alpha > 0$, thus

$$\mathfrak{F}_T = (\mathfrak{F}_{T^\alpha})^{1/\alpha} \stackrel{\text{def}}{=} \{\mathcal{L}(W): \mathcal{L}(W^\alpha) \in \mathfrak{F}_{T^\alpha}\} \tag{1.7}$$

for any $\alpha > 0$, where $T^\alpha \stackrel{\text{def}}{=} (T_j^\alpha)_{j \in J}$. As in [4], the unique $-\alpha$ solving (1.5) will be called the *characteristic exponent* of $T$. Note that $\mathcal{L}(W) = \text{Weib}(c, \alpha)$, the Weibull distribution with distribution function $(1 - \exp(-ct^\alpha))\mathbf{1}_{(0,\infty)}(t)$, if $\mathcal{L}(W^\alpha) = \text{Exp}(c)$.

Equation (1.1) with positive $T_j$ may be rewritten in terms of the distribution function $F(t) \stackrel{\text{def}}{=} \mathbb{P}(W < t)$ as

$$1 - F(t) = \mathbb{P}(T_j W_j \geq t \text{ for all } j \in J) = \prod_{j \in J}\left(1 - F\left(\frac{t}{T_j}\right)\right), \quad t \in \mathbb{R}. \tag{1.8}$$

Choosing $t = 0$ and using $|J| \geq 2$, we arrive at the basic conclusion that

$$1 - F(0) \leq (1 - F(0))^2$$

that is $F(0) \in \{0, 1\}$. Since furthermore

$$1 - F(t+) \leq \mathbb{P}(T_j W_j > t \text{ for all } j \in J) = \prod_{j \in J}\left(1 - F\left(\left(\frac{t}{T_j}\right)+\right)\right), \quad t \in \mathbb{R}, \tag{1.9}$$

we also infer $F(0+) \in \{0, 1\}$ and thus $F(0+) - F(0) \in \{0, 1\}$. Any fixed point $F \neq \delta_0$ is hence necessarily continuous at 0 and concentrated on either $\mathbb{R}^>$ (positive solution) or $\mathbb{R}^< \stackrel{\text{def}}{=} (-\infty, 0)$ (negative solution). Denote by $\mathfrak{F}_T^+$ and $\mathfrak{F}_T^-$ the set of, respectively, positive and negative solutions to (1.1) and notice that the trivial fixed point $\delta_0$ is contained in neither of these two sets. An application of the idempotent transformation $t \mapsto -\frac{1}{t}$ from $\mathbb{R}^* \stackrel{\text{def}}{=} \mathbb{R} \setminus \{0\}$ to $\mathbb{R}^*$ immediately shows that $\mathcal{L}(W) \neq \delta_0$ is a solution to (1.1) iff $\mathcal{L}(-1/W)$ is a solution to the very same equation with weight vector $T^{-1}$ instead of $T$, i.e.

$$\mathfrak{F}_T \setminus \{\delta_0\} = \left\{\mathcal{L}\left(-\frac{1}{W}\right): \mathcal{L}(W) \in \mathfrak{F}_{T^{-1}} \setminus \{\delta_0\}\right\}. \tag{1.10}$$

Moreover, a negative solution $F \in \mathfrak{F}_T$ always corresponds to a positive one in $\mathfrak{F}_{T^{-1}}$, so

$$\mathfrak{F}_T^- = \left\{\mathcal{L}\left(-\frac{1}{W}\right): \mathcal{L}(W) \in \mathfrak{F}_{T^{-1}}^+\right\}. \tag{1.11}$$

It therefore suffices to determine the positive solutions to (1.1) for arbitrary positive weight vector $T$.

The further organization of this article is as follows. Section 2 is mainly devoted to a discussion of all those cases where $\mathfrak{F}_T^+$ is fairly simply identified. At the end of Section 2 we will be left with only two cases, namely when $\inf_{j \in J} T_j > 1$ and

$$\text{either } |J| < \infty \quad \text{or} \quad J = \mathbb{N}, \lim_{j \to \infty} T_j = \infty, \tag{1.12}$$

which require indeed a deeper analysis given in Section 3. Its main result, Proposition 3.3, provides us with the crucial information that under (1.12) the characteristic exponent of $T$ necessarily exists and is negative unless $\mathfrak{F}_T^+$ is empty. The analysis further shows that any $F \in \mathfrak{F}_T^+$ must have unbounded support which indeed forms its most intricate part. We thus conclude that $\nu_F(t) \stackrel{\text{def}}{=} -\log(1 - F(t))$ defines a Radon measure on $\mathbb{R}^>$ which, by (1.8), satisfies

$$\nu_F(t) = \sum_{j \in J} \nu_F\left(\frac{t}{T_j}\right) = \nu_F \star \Lambda(t), \quad t \in \mathbb{R}^>, \tag{1.13}$$

and is hence $\Lambda$-harmonic, where $\Lambda \stackrel{\text{def}}{=} \sum_{j \in J} \delta_{T_j}$ and $\star$ means multiplicative convolution. This crucial observation will bring Choquet theory into play in a very similar manner as in [4] for another type of stochastic fixed point equation, namely

$$W \stackrel{\text{d}}{=} \sum_{j \in J} T_j W_j \tag{1.14}$$

for a sequence $T = (T_j)_{j \in J}$ of nonzero real-valued numbers. The powerful Choquet–Deny theorem [5] will enable us to obtain very explicit information on the possible form of $\nu_F$ and thus $F$ which is subsequently converted into a full description of $\mathfrak{F}_T^+$ which can be found in Section 4. Having solved the case (C1) completely, the final section provides the results for the simpler cases (C2) and (C3).

## 2. Basic results and simple cases

As already announced in the Introduction, the present section will collect some basic results about solutions to (1.1) including a discussion of those cases where $\mathfrak{F}_T^+$, or even $\mathfrak{F}_T$, is easily determined. Of course, $|J| \geq 2$ will always be in force throughout without further mention.

Given a distribution $F$ on $\mathbb{R}^>$, we always use the same symbol for its left continuous distribution function and put further $\overline{F}(t) \stackrel{\text{def}}{=} 1 - F(t)$ and

$$\nu_F(t) \stackrel{\text{def}}{=} -\log \overline{F}(t) \tag{2.1}$$

for $t \in \mathbb{R}$. Notice that $\nu_F(t)$ is finite, nondecreasing and left continuous on $[0, u_F)$ and positive on $(l_F, u_F)$, where $l_F \stackrel{\text{def}}{=} \sup\{t \geq 0 \colon F(t) = 0\}$ and $u_F \stackrel{\text{def}}{=} \inf\{t \geq 0 \colon F(t) = 1\}$. Moreover, $\lim_{t \downarrow l_F} \nu_F(t) = 0$ and $\lim_{t \uparrow u_F} \nu_F(t) = \infty$. Consequently, $\nu_F$ defines a Radon measure on $(0, u_F)$.

Defining $L(v) \stackrel{\text{def}}{=} \prod_{j=1}^n T_{v_j}$ for $v = (v_1, \ldots, v_n) \in J^n$ and $n \geq 1$, an iteration of (1.1) leads to

$$W \stackrel{\text{d}}{=} \min_{v \in J^n} L(v) W(v), \tag{2.2}$$

for each $n \geq 1$, where the $W(v)$ are i.i.d. copies of $W$. This is the weighted branching representation of $\mathcal{L}(W)$ because the $L(v)$ may be interpreted as the total weight of the branch $\emptyset \to (v_1) \to (v_1, v_2) \to \cdots \to v = (v_1, \ldots, v_n)$ in the Ulam–Harris tree $\bigcup_{n \geq 0} J^n$, where $J^0 \stackrel{\text{def}}{=} \{\emptyset\}$. An edge from $v$ to $(v, j)$, $j \in J$, carries the weight $T_j$, and the total weight of a branch is obtained by multiplication. The corresponding equation for the distribution function $F$ of $W$ takes the form (compare (1.8))

$$\overline{F}(t) = \prod_{v \in J^n} \overline{F}\left(\frac{t}{L(v)}\right), \quad t \in \mathbb{R}. \tag{2.3}$$

Our first lemma shows that $\inf_{j \in J} T_j \geq 1$ forms a necessary condition for the existence of positive solutions to (1.1).

**Lemma 2.1.** *If $\inf_{j \in J} T_j < 1$ then $\mathfrak{F}_T^+ = \emptyset$.*

**Proof.** W.l.o.g. suppose $T_1 < 1$. Given any $F \in \mathfrak{F}_T$, (2.3) implies

$$\overline{F}(t) = \prod_{v \in J^n} \overline{F}\left(\frac{t}{L(v)}\right) \leq \overline{F}\left(\frac{t}{T_1^n}\right)$$

for all $n \geq 1$ and $t > 0$, whence

$$\overline{F}(t) \leq \lim_{n \to \infty} \overline{F}\left(\frac{t}{T_1^n}\right) = 0$$

for all $t > 0$, i.e. $F(0+) = 1$. □

By combining Lemma 2.1 with (1.11), we see that $\sup_{j \in J} T_j > 1$ implies $\mathfrak{F}_T^- = \emptyset$, and with this observation the following corollary is immediate.

**Corollary 2.2.** *Suppose $\inf_{j \in J} T_j \geq 1$. Then $\mathfrak{F}_T = \{\delta_c \colon c \in \mathbb{R}\}$, if $T_j = 1$ for all $j \in J$, while $\mathfrak{F}_T = \{\delta_0\} \cup \mathfrak{F}_T^+$, otherwise.*

In view of the previous results it is clear that positive solutions to (1.1) can only occur if $\inf_{j \in J} T_j \geq 1$ which is therefore assumed hereafter unless stated otherwise. A further analysis requires the distinction of several subcases listed as (A1)–(A6) below:

(A1) $|\{j \in J \colon T_j = 1\}| \geq 2$.
(A2) $J = \mathbb{N}$, $|\{j \in J \colon T_j = 1\}| \leq 1$ and $\liminf_{j \to \infty} T_j = 1$.
(A3) $|\{j \in J \colon T_j = 1\}| = 1$ and $\inf_{j \in J^*} T_j > 1$, where $J^* \stackrel{\text{def}}{=} \{j \in J \colon T_j \neq 1\}$.

(A4) $J = \mathbb{N}$, $\inf_{j \in J} T_j > 1$ and $\liminf_{j \to \infty} T_j < \infty$.
(A5) $|J| < \infty$ and $\inf_{j \in J} T_j > 1$.
(A6) $J = \mathbb{N}$, $\inf_{j \in J} T_j > 1$ and $\lim_{j \to \infty} T_j = \infty$.

Plainly, (A1)–(A3) are subcases of $\inf_{j \in J} T_j = 1$, while (A4)–(A6) are subcases of $\inf_{j \in J} T_j > 1$. As for (A1)–(A4), a complete description of $\mathfrak{F}_T^+$ is rather easily obtained and stated in the following proposition. The characteristic exponent of $T$ and the use of Choquet theory enter for the remaining cases (A5) and (A6), the last case being the most difficult one as involving harmonic analysis on trees. The latter will be provided in Section 3 followed by the description of $\mathfrak{F}_T^+$ in Section 4.

**Proposition 2.3.** *Suppose $\inf_{j \in J} T_j \geq 1$.*

(a) *If (A1) or (A2) holds true, then $\mathfrak{F}_T^+ = \{\delta_c \colon c > 0\}$.*
(b) *If (A3) holds true, then $\mathfrak{F}_T^+$ consists of all distributions $F$ with $0 < l_F \leq u_F < \infty$ such that $\inf_{j \in J^*} T_j \geq u_F / l_F$ (and thus includes $\{\delta_c \colon c > 0\}$).*
(c) *If (A4) holds true, then $\mathfrak{F}_T^+$ is empty.*

**Proof.** (a) Given (A1), we infer from (1.8) for any $F \in \mathfrak{F}_T^+$ that

$$\overline{F}(t) \leq \overline{F}(t)^2$$

and thus $F(t) \in \{0, 1\}$ for all $t > 0$, that is $F = \delta_c$ for some $c > 0$. If (A2) holds true then $F \in \mathfrak{F}_T^+$ still satisfies the slightly weaker inequality

$$\overline{F}(t) \leq \overline{F}(t - \varepsilon)^2$$

for all $t > 0$ and $\varepsilon > 0$ which leads to the same conclusion upon letting $\varepsilon$ tend to 0 and using the left continuity of $F$.

(b) W.l.o.g. let $T_1 = 1$ and put $c \stackrel{\text{def}}{=} \inf_{j \geq 2} T_j$. If $F$ is any distribution with $0 < l_F \leq u_F < \infty$ such that $c \geq u_F / l_F$ and if $W_1, W_2, \ldots$ are i.i.d. with distribution $F$ then we obviously have

$$\inf_{j \geq 2} T_j W_j \geq c l_F \geq u_F \geq W_1 = T_1 W_1 \quad \text{a.s.}$$

and thus

$$\inf_{j \geq 1} T_j W_j \stackrel{\text{d}}{=} W_1,$$

i.e. $F \in \mathfrak{F}_T^+$. Conversely, given any positive fixed point $F$, Eq. (1.8) yields for each $0 < t < u_F$ ($\Rightarrow \overline{F}(t) > 0$)

$$\overline{F}(t) = \overline{F}(t) \prod_{j \geq 2} \overline{F}\left(\frac{t}{T_j}\right)$$

and thereby (recalling the left continuity of $F$)

$$1 = \prod_{j \geq 2} \overline{F}\left(\frac{t}{T_j}\right) \leq \inf_{j \geq 2} \overline{F}\left(\frac{t}{T_j}\right) = \overline{F}\left(\frac{t}{c}\right),$$

hence $0 < u_F < \infty$ and $l_F \geq u_F / c > 0$.

(c) Suppose there is a positive fixed point $F$. Put $c_1 \stackrel{\text{def}}{=} \inf_{j \geq 1} T_j$, $c_2 \stackrel{\text{def}}{=} \liminf_{j \to \infty} T_j$ (finite by assumption) and $J_\varepsilon \stackrel{\text{def}}{=} \{j \geq 1 \colon T_j \geq c_2 - \varepsilon\}$ for any $\varepsilon \in (0, c_2 - 1)$. As $F \neq \delta_0$ there must be a $t > 0$ with $\overline{F}(t) > 0$. Using $|J_\varepsilon| = \infty$ we infer from (1.8)

$$0 < \overline{F}(t) \leq \prod_{j \in J_\varepsilon} \overline{F}\left(\frac{t}{T_j}\right) \leq \overline{F}\left(\frac{t}{c_2 - \varepsilon}\right)^n$$

for every $n \geq 1$ and thus $\overline{F}(\frac{t}{c_2-\varepsilon}) = 1$, which in turn entails $l_F > 0$. On the other hand we may then pick $\eta > 0$ such that $\frac{l_F+\eta}{c_1} \leq l_F$ and are thus led to the contradiction

$$1 > \overline{F}(l_F + \eta) = \prod_{j \geq 1} \overline{F}\left(\frac{l_F+\eta}{T_j}\right) \geq \prod_{j \geq 1} \overline{F}\left(\frac{l_F+\eta}{c_1}\right) = 1.$$

Hence a positive fixed point cannot exist. □

Let us finally mention that in each of the four considered cases (A1)–(A4) the characteristic exponent of $T$ does obviously not exist.

## 3. The characteristic exponent of $T$

The present section deals with the most interesting cases (A5) and (A6) and provides the necessary results needed to determine $\mathfrak{F}_T^+$ to be done in Section 4. Recall from the Introduction that $\Lambda = \sum_{j \in J} \delta_{T_j}$, $\nu_F(t) = -\log \overline{F}(t)$ and that, by (1.13), $\nu_F$ is $\Lambda$-harmonic. Further note that $\inf_{j \in J} T_j > 1$ implies that any positive solution is necessarily nontrivial. The following lemma further shows that it must carry mass in any neighborhood of 0.

**Lemma 3.1.** *If* (A5) *or* (A6) *holds then any* $F \in \mathfrak{F}_T^+$ *satisfies* $l_F = 0$. *Moreover,* (A5) *also implies* $u_F = \infty$.

**Proof.** Put $c \stackrel{\text{def}}{=} \inf_{j \in J} T_j$. Equation (1.1) with $W, W_1, W_2, \ldots \stackrel{\text{d}}{=} F$, the definition of $l_F$ and

$$\mathbb{P}(W \geq cl_F) = \mathbb{P}\left(\inf_{j \in J} T_j W_j \geq cl_F\right) = 1$$

together imply $l_F \geq cl_F$ and thus $l_F = 0$ because $c > 1$.

If (A5) and thus $|J| < \infty$ holds true then furthermore, by (2.3),

$$\overline{F}(t) = \prod_{v \in J^n} \overline{F}\left(\frac{t}{L(v)}\right) \geq \overline{F}\left(\frac{t}{c^n}\right)^{|J|^n}$$

for all $t > 0$ and $n \geq 1$. As $l_F = 0$ we can choose $n = n(t)$ so large that $\overline{F}(t/c^n) > 0$ and thus $\overline{F}(t) > 0$. This proves $u_F = \infty$. □

That under (A6) any $F \in \mathfrak{F}_T^+$ satisfies $u_F = \infty$, too, is more difficult to verify and in fact derived in Proposition 3.3 below proved by harmonic analysis. The importance of $u_F = \infty$ stems from the fact that only then the $\Lambda$-harmonic $\nu_F$ defines a Radon measure on $\mathbb{R}^>$ which in turn forms a crucial requirement for the use of Choquet theory needed to identify the form of $\nu_F$ and thus $F$.

Define the function $m : \mathbb{R} \to (0, \infty]$ by

$$m(\beta) \stackrel{\text{def}}{=} \sum_{j \in J} T_j^\beta. \tag{3.1}$$

Note that $m$ is continuous and, as all $T_j$ are $> 1$, strictly increasing on $\{\beta : m(\beta) < \infty\}$ with $m(0) = |J| \geq 2$. Hence the characteristic exponent of $T$, if it exists, is necessarily unique and negative because $m(0) > 1$. If (A5) holds, $m$ is everywhere finite with $\lim_{\beta \to -\infty} m(\beta) = 0$ and $\lim_{\beta \to \infty} m(\beta) = \infty$. Consequently, the characteristic exponent of $T$ exists. Turning to (A6), the latter may fail but Proposition 3.3 shows that then $\mathfrak{F}_T^+$ is empty.

With view to the subsequent applications of the Choquet–Deny theorem [5] we continue with the collection of some necessary facts. Let $\mathbb{G}_\Lambda$ denote the closed multiplicative subgroup of $\mathbb{R}^>$ generated by $\Lambda$. Our standing assumption $|J| \geq 2$ excludes the trivial subgroup $\{1\}$ so that either

$$\mathbb{G}_\Lambda = \mathbb{R}^> \quad \text{(continuous case)},$$

or

$$\mathbb{G}_\Lambda = r^\mathbb{Z} \quad \text{for some } r > 1, \text{ where } r^\mathbb{Z} \stackrel{\text{def}}{=} \{r^z \colon z \in \mathbb{Z}\} \quad (r\text{-geometric case}).$$

The Haar measure (unique up to multiplicative constants), denoted as $\lambda_{\mathbb{G}_\Lambda}$ hereafter, equals $|u|^{-1}\,du$ in the continuous case and counting measure in the $r$-geometric one. Let $E(\mathbb{G}_\Lambda)$ be the set of characters of $\mathbb{G}_\Lambda$, that is the set of all continuous positive functions $e \colon \mathbb{G}_\Lambda \to \mathbb{R}^>$ satisfying $e(xy) = e(x)e(y)$ for all $x, y \in \mathbb{G}_\Lambda$. Of particular interest for our purposes is the subset

$$E_1(\Lambda) \stackrel{\text{def}}{=} \left\{ e \in E(\mathbb{G}_\Lambda) \colon \int e(x^{-1}) \Lambda(dx) = 1 \right\}.$$

It is not difficult to check that in both cases the characters are given by the functions $e_\beta(x) \stackrel{\text{def}}{=} |x|^{-\beta}$, $\beta \in \mathbb{R}$, so $E(\mathbb{G}_\Lambda) = E$ is independent of $\Lambda$. Moreover, we see upon noting $\int e_\beta(x^{-1})\Lambda(dx) = \sum_{j \geq 1} T_j^\beta = m(\beta)$ that $E_1(\Lambda)$ is either void or consists of the single element $e_{-\alpha}$, $-\alpha < 0$ the characteristic exponent of $T$. Hence, $E_1(\Lambda) = \{e_{-\alpha}\}$ always holds true in the case (A5), whereas $E_1(\Lambda) = \emptyset$ may happen under (A6).

Now consider a Radon measure $\mu$ on $\mathbb{R}^>$ and suppose that $\mu$ is $\Lambda$-*harmonic*, defined by $\mu = \mu \star \Lambda$. Here $\star$ means multiplicative convolution, that is

$$\int f(x) \mu \star \Lambda(dx) \stackrel{\text{def}}{=} \iint f(xy) \mu(dx) \Lambda(dy)$$

for any measurable $f \colon \mathbb{G}_\Lambda \to [0, \infty)$. The set of all $\Lambda$-harmonic measures is a convex cone. By the Choquet–Deny theorem [5] we infer that any nonzero $\Lambda$-harmonic $\mu$ has a unique integral representation

$$\mu = \int \mu_e(y^{-1} \cdot) \overline{\mu}(de, dy),$$

where $\mu_e(dx) \stackrel{\text{def}}{=} e(x) \lambda_{\mathbb{G}_\Lambda}(dx)$ for $e \in E$ and $\overline{\mu}$ is a finite measure on $E_1(\Lambda) \times \mathbb{R}^>/\mathbb{G}_\Lambda$ endowed with the Baire $\sigma$-field. If $E_1(\Lambda) = \emptyset$ there is no $\Lambda$-harmonic measure. Otherwise, $E_1(\Lambda) = \{e_{-\alpha}\}$ so that $\overline{\mu}$ must equal $c(\delta_{e_{-\alpha}} \otimes \mu)$ for some probability measure $\mu$ on the factor group $\mathbb{R}^>/\mathbb{G}_\Lambda$ and a $c > 0$. This means that

$$\mu(dx) = \int_{\mathbb{R}^>/\mathbb{G}_\Lambda} c \left(\frac{x}{y}\right)^\alpha \lambda_{y\mathbb{G}_\Lambda}(dx) \mu(dy). \tag{3.2}$$

Of course, if $\mathbb{G}_\Lambda = \mathbb{R}^>$ then $\mathbb{R}^>/\mathbb{G}_\Lambda = \{1\}$, $\mu = \delta_1$ and thus

$$\mu(dx) = c x^{\alpha - 1} \mathbf{1}_{\mathbb{R}^>}(x)\, dx \tag{3.3}$$

for some $c > 0$. In the discrete case where $\mathbb{G}_\Lambda = r^\mathbb{Z}$ and $\mathbb{R}^>/\mathbb{G}_\Lambda = [1, r)$ for some $r > 1$ it follows that

$$\mu(dx) = c \int_{[1,r)} \sum_{n \in \mathbb{Z}} r^{n\alpha} \delta_{yr^n}(dx) \mu(dy) \tag{3.4}$$

for some $c > 0$ and a probability measure $\mu$ on $[1, r)$.

The next lemma is stated for later reference and provides us with the general form of $\nu_F$ in the simpler case (A5).

**Lemma 3.2.** *If* (A5) *holds, then the characteristic exponent* $-\alpha < 0$ *of* $T$ *exists and* $\nu_F$ *is of the form* (3.3) *or* (3.4), *i.e.*

$$\nu_F(\mathrm{d}x) = \begin{cases} c\mathbf{1}_{(0,\infty)}(x)\, x^{\alpha-1}\,\mathrm{d}x & \text{if } \mathbb{G}_\Lambda = \mathbb{R}^>, \\ \int_{[1,r)} \sum_{n\in\mathbb{Z}} c r^{n\alpha} \delta_{yr^n}(\mathrm{d}x)\hat{F}(\mathrm{d}y) & \text{if } \mathbb{G}_\Lambda = r^{\mathbb{Z}} \text{ for } r > 1, \end{cases} \tag{3.5}$$

*for some* $c > 0$ *(and a probability measure* $\hat{F}$ *on* $[1,r)$ *in the discrete case).*

**Proof.** We already noted above that $T$ possesses a unique negative characteristic exponent $-\alpha$ because $m(0) = |J| \geq 2$ and $J$ is finite. Hence (3.5) is a direct consequence of (3.3) and (3.4) recalling once more that $\nu_F$ is a $\Lambda$-harmonic Radon measure. □

**Proposition 3.3.** *If* (A6) *holds and* $\mathfrak{F}_T^+$ *is not empty, then the characteristic exponent of* $T$ *exists and any* $F \in \mathfrak{F}_T^+$ *satisfies* $u_F = \infty$.

**Proof.** Since, by assumption, $\lim_{j\to\infty} T_j = \infty$ it is no loss of generality to assume that $1 < T_1 \leq T_2 \leq \cdots$. Let $F \in \mathfrak{F}_T^+$. Then $l_F = 0$, by Lemma 3.1, and $u_F > 0$. For $t \in (0, u_F)$ and $\lambda \geq 0$ put

$$g(\lambda, t) \stackrel{\text{def}}{=} \frac{\log \overline{F}(\lambda t)}{\log \overline{F}(t)} = \frac{\nu_F(\lambda t)}{\nu_F(t)} \tag{3.6}$$

and note that $g(\lambda, t) = \infty$ if $u_F < \infty$ and $\lambda \geq u_F/t$. To show that $u_F = \infty$ and thus $g(\lambda, t)$ is always finite will be one of the difficult tasks of this proof. By definition of $l_F$ and $u_F$, for each $t \in (0, u_F)$:

(1) $g(0, t) = 0$, $g(1, t) = 1$ and $0 < g(\lambda, t) \leq 1$ for $0 < \lambda < 1$.
(2) $g(\cdot, t)$ is left continuous and nondecreasing on $[0, u_F/t)$.

The fixed point Eq. (1.8) (or (1.13)) implies for all $t \in (0, u_F)$ and $\lambda \geq 0$

$$g(\lambda, t) = \sum_{j \geq 1} g\!\left(\frac{\lambda}{T_j}, t\right) = \sum_{j \geq 1} g\!\left(\lambda, \frac{t}{T_j}\right) g\!\left(\frac{1}{T_j}, t\right) \tag{3.7}$$

as well as

$$\sum_{j \geq 1} g\!\left(\frac{1}{T_j}, t\right) = 1. \tag{3.8}$$

More generally, we have upon iteration and using the weighted branching representation described in Section 2,

$$g(\lambda, t) = \sum_{|v|=n} g\!\left(\frac{\lambda}{L(v)}, t\right) = \sum_{|v|=n} g\!\left(\lambda, \frac{t}{L(v)}\right) g\!\left(\frac{1}{L(v)}, t\right) \tag{3.9}$$

and

$$\sum_{|v|=n} g\!\left(\frac{1}{L(v)}, t\right) = 1 \tag{3.10}$$

for every $n \geq 1$, where $|v|$ denotes the length of the vector $v \in \mathbb{V} \stackrel{\text{def}}{=} \{\emptyset\} \cup \bigcup_{n \geq 1} \mathbb{N}^n$. We further use $vj$ and $vw$ as shorthand notation for $(v, j)$ and $(v, w)$ if $j \geq 1$ and $v, w \in \mathbb{V}$.

Fix any $t \in (0, u_F)$ and let $M = (M_n)_{n \geq 0}$ be a Markov chain on a probability space $(\Omega, \mathfrak{U}, \mathbb{P})$ with state space $\mathbb{V}$ and 1-step transition kernel $P = P_t$ defined by

$$P(v, \{vj\}) = g\!\left(\frac{1}{T_j}, \frac{t}{L(v)}\right), \quad j \geq 1, v \in \mathbb{V}. \tag{3.11}$$

One can easily check that the $n$-step transition kernel $P^n$ satisfies

$$P^n(v, \{vw\}) = g\left(\frac{1}{L(w)}, \frac{t}{L(v)}\right), \quad j \geq 1, v, w \in \mathbb{V}, |w| = n. \tag{3.12}$$

Let $\mathbb{P}$ be such that $\mathbb{P}(M_0 = v) > 0$ for all $v \in \mathbb{V}$ and put $\mathbb{P}_v \stackrel{\text{def}}{=} \mathbb{P}(\cdot | M_0 = v)$. Since, for any $v \in \mathbb{V}$ and $\lambda \in [0, u_F L(v)/t)$,

$$\mathbb{E}_v g\left(\lambda, \frac{t}{L(M_1)}\right) = \sum_{j \geq 1} g\left(\lambda, \frac{t}{L(v)T_j}\right) g\left(\frac{1}{T_j}, \frac{t}{L(v)}\right) = g\left(\lambda, \frac{t}{L(v)}\right), \tag{3.13}$$

and since $L(w) \geq T_1^{|w|} \to \infty$ as $|w| \to \infty$, we see that $(g(\lambda, t/L(M_n)))_{n \geq n(\lambda)}$ forms a (bounded) $\mathbb{P}_v$-martingale for each $v \in \mathbb{V}$ and $\lambda \geq 0$, where $n(\lambda)$ is chosen so large that $\lambda < u_F T_1^{n(\lambda)}/t \leq u_F L(M_{n(\lambda)})/t$. Denote by $Y(\lambda)$ its a.s. and $L^1$-limit under $\mathbb{P}$. Notice that

$$\mathbb{E}_v Y(\lambda) = g\left(\lambda, \frac{t}{L(v)}\right) \tag{3.14}$$

for every $v \in \mathbb{V}$ and $\lambda \in [0, u_F L(v)/t)$. By using the monotonicity and left continuity of $g(\cdot, t)$, it is not difficult to verify that outside a $\mathbb{P}$-null set $\mathcal{N}$

$$\lim_{n \to \infty} g\left(\lambda, \frac{t}{L(M_n)}\right) = Y(\lambda) \quad \mathbb{P}\text{-a.s.}$$

holds simultaneously *for all* $\lambda \geq 0$, and $Y$ is a finite, nondecreasing, left continuous random function on $\mathcal{N}^c$ satisfying $Y(0) = 0$ and $Y(1) = 1$. We claim that this in combination with (3.7) and (3.14) implies

$$Y(\lambda) = \sum_{j \geq 1} Y\left(\frac{\lambda}{T_j}\right) \quad \mathbb{P}\text{-a.s.} \tag{3.15}$$

for all $\lambda \geq 0$. Indeed, we have by (3.7) and Fatou's lemma

$$\sum_{j \geq 1} Y\left(\frac{\lambda}{T_j}\right) \leq \lim_{n \to \infty} \sum_{j \geq 1} g\left(\frac{\lambda}{T_j}, \frac{t}{L(M_n)}\right) = \lim_{n \to \infty} g\left(\lambda, \frac{t}{L(M_n)}\right) = Y(\lambda) \quad \mathbb{P}\text{-a.s.}$$

(since $\frac{t}{L(M_n)} \in (0, u_F)$ and $\lambda \in [0, \frac{u_F L(M_n)}{t})$ for all $n$ large enough), while (3.14) and (3.7) give

$$\sum_{j \geq 1} \mathbb{E}_v Y\left(\frac{\lambda}{T_j}\right) = \sum_{j \geq 1} g\left(\frac{\lambda}{T_j}, \frac{t}{L(v)}\right) = g\left(\lambda, \frac{t}{L(v)}\right) = \mathbb{E}_v Y(\lambda)$$

for all $(v, \lambda) \in \mathbb{V} \times [0, \infty)$ having $\lambda < L(v)/t$. We infer

$$\mathbb{P}_v\left(Y(\lambda) = \sum_{j \geq 1} Y\left(\frac{\lambda}{T_j}\right)\right) = 1$$

for all such $v$ and $\lambda$. In order to obtain the very same for all $(v, \lambda) \in \mathbb{V} \times [0, \infty)$ (i.e. (3.15)), use that $L(v) \to \infty$ as $|v| \to \infty$ and that, by (3.12), the càglàd process $Y = (Y(\lambda)_{\lambda \geq 0}$ satisfies

$$\mathbb{P}_v(Y \in \cdot) = \sum_{|w|=n} P^n(v, \{vw\}) \mathbb{P}_{vw}(Y \in \cdot) = \sum_{|w|=n} g\left(\frac{1}{L(w)}, \frac{t}{L(v)}\right) \mathbb{P}_{vw}(Y \in \cdot),$$

where $g(\frac{1}{L(w)}, \frac{t}{L(v)}) > 0$ for all $v, w \in \mathbb{V}$. Consequently, fixing any $v \in \mathbb{V}$ and $\lambda \geq 0$, we conclude upon choosing $n$ so large that $\lambda < L(vw)/t$ for all $w$ with $|w| \geq n$

$$\mathbb{P}_v\left(Y(\lambda) = \sum_{j \geq 1} Y\left(\frac{\lambda}{T_j}\right)\right) = \sum_{|w|=n} g\left(\frac{1}{L(w)}, \frac{t}{L(v)}\right) \mathbb{P}_{vw}\left(Y(\lambda) = \sum_{j \geq 1} Y\left(\frac{\lambda}{T_j}\right)\right) = 1.$$

We must still prove that $Y(\lambda) > 0$ $\mathbb{P}$-a.s. for all $\lambda > 0$. But

$$1 = Y(1) = \sum_{j \geq 1} Y\left(\frac{1}{T_j}\right) \quad \text{a.s.}$$

together with $Y(1/T_1) = \sup_{j \geq 1} Y(1/T_j)$ $\mathbb{P}$-a.s. implies $Y(1/T_1) > 0$ and thus $Y(\lambda) > 0$ $\mathbb{P}$-a.s. for all $\lambda \geq 1/T_1$. Repeating the above argument with $Y(1/T_1)$ instead of $Y(1)$ gives $Y(1/T_1^2) > 0$ and thus $Y(\lambda) > 0$ $\mathbb{P}$-a.s. for all $\lambda \geq 1/T_1^2$. Continuing this way the assertion easily follows.

The main conclusion from the previous analysis is that for $\mathbb{P}^M$ almost every infinite path $x = v_1 v_2 \cdots \in \partial \mathbb{V}$ we have upon setting $L_n(x) \stackrel{\text{def}}{=} L(v_1 \cdots v_n)$ for $n \geq 1$ and $L_0(x) \stackrel{\text{def}}{=} 1$ that

$$G(\lambda, x) \stackrel{\text{def}}{=} \lim_{n \to \infty} g\left(\lambda, \frac{t}{L_n(x)}\right) \tag{3.16}$$

exists for all $\lambda \geq 0$, is nondecreasing, left continuous with $G(0, x) = 0$, $G(1, x) = 1$ and $0 < G(\lambda, x) < \infty$ for all $\lambda > 0$. Moreover, it satisfies

$$G(\lambda, x) = \sum_{j \geq 1} G\left(\frac{\lambda}{T_j}, x\right) = G(\cdot, x) \star \Lambda(\lambda)$$

for all $\lambda \geq 0$ which means that $G(\cdot, x)$, viewed as a Radon measure on $\mathbb{R}^>$, is $\Lambda$-harmonic. Hence an application of the Choquet–Deny theorem ensures that the characteristic exponent $-\alpha$ of $T$ exists and that $G(\cdot, x)$ is of the form described in (3.3) or (3.4). This means that

$$G(\lambda, x) = c(x) \lambda^\alpha,$$

if $\mathbb{G}_\Lambda = \mathbb{R}^>$, and

$$G(\lambda, x) = h(\lambda, x) \lambda^\alpha,$$

if $\mathbb{G}_\Lambda = r^{\mathbb{Z}}$ for some $r > 1$. Here $h(\cdot, x) : \mathbb{R}^> \to [0, \infty)$ denotes a multiplicatively $r$-periodic, left continuous function such that $h(\lambda, x) \lambda^\alpha$ is nondecreasing in $\lambda$ (see at the beginning of Section 4 for details). The normalization $G(1, x) = 1$ implies $c(x) \equiv 1$ if $\mathbb{G}_\Lambda = \mathbb{R}^>$ and $h(1, x) \equiv 1$ if $\mathbb{G}_\Lambda = r^{\mathbb{Z}}$. Use the $L^1$-convergence of $g(\lambda, t/L(M_n))$ to $Y(\lambda)$ for each $0 \leq \lambda < u_F/t$ to infer

$$g(\lambda, t) = \lim_{n \to \infty} \mathbb{E}_\emptyset g\left(\lambda, \frac{t}{L(M_n)}\right) = \mathbb{E}_\emptyset Y(\lambda) = \begin{cases} \lambda^\alpha & \text{if } \mathbb{G}_\Lambda = \mathbb{R}^>, \\ \overline{h}(\lambda) \lambda^\alpha & \text{if } \mathbb{G}_\Lambda = r^{\mathbb{Z}}, \end{cases} \tag{3.17}$$

where $\overline{h}(\lambda) \stackrel{\text{def}}{=} \mathbb{E}_\emptyset h(\lambda, t/L(M_\infty))$ and $M_\infty \stackrel{\text{def}}{=} \lim_{n \to \infty} M_n$ in the usual topology on $\mathbb{V} \cup \partial \mathbb{V}$. It is now an easy step to show that (3.17) actually holds for all $\lambda \geq 0$. We do so only for the case $\mathbb{G}_\Lambda = \mathbb{R}^>$ because the modifications in the discrete case are very similar. By the first equality in (3.6),

$$g(\lambda T_1, t) = \sum_{j \geq 1} g\left(\frac{\lambda T_1}{T_j}, t\right) = \sum_{j \geq 1} \left(\frac{\lambda T_1}{T_j}\right)^\alpha = (\lambda T_1)^\alpha \sum_{j \geq 1} T_j^{-\alpha} = (\lambda T_1)^\alpha,$$

for all $0 \leq \lambda < u_F/t$, where we have used that $T_1 = \min_j T_j$ and $\sum_j T_j^{-\alpha} = 1$. So $g(\lambda, t)$ satisfies (3.17) for all $0 \leq \lambda < u_F T_1/t$. An inductive argument yields (3.17) for all $0 \leq \lambda < u_F T_1^n/t$ and $n \geq 1$ and thus for all $\lambda \geq 0$. A particular consequence is that $\nu_F(\lambda t) < \infty$ for all $\lambda \geq 0$ and therefore $u_F = \infty$. The proof is herewith complete. □

We finally note as a consequence of the previous proof that $\nu_F$ is again of the form (3.5) in the case (A6).

## 4. Case (C1): Results in the nontrivial subcases

In order to state our result for the remaining cases (A5) and (A6) we must briefly introduce a class of distributions which turn up there in cases where $\mathbb{G}_\Lambda = r^{\mathbb{Z}}$ for some $r > 1$. A function $h : \mathbb{R}^{>} \to [0, \infty)$ is called *multiplicatively r-periodic* for $r > 1$ if $h(rt) = h(t)$ for all $t > 0$. Let $\mathfrak{H}(r, \alpha)$ be the class of all left continuous functions of this type such that $h(t)t^\alpha$ is nondecreasing.

**Definition 4.1.** *Let $r > 1$, $\alpha > 0$ and $h \in \mathfrak{H}(r, \alpha)$. The distribution $F$ on $\mathbb{R}^{>}$ with survival function*

$$\overline{F}(t) = \exp(-h(t)t^\alpha), \quad t > 0, \tag{4.1}$$

*is called r-periodic Weibull distribution with parameters $h$ and $\alpha$, in short r-Weib$(h, \alpha)$. In case $\alpha = 1$ it is also called r-periodic exponential distribution with parameter $h$, in short r-Exp$(h)$.*

These distributions may be viewed as periodic variants of their familiar continuous counterparts on $\mathbb{R}^{>}$. The particular choice $h \in \mathfrak{H}(r, \alpha)$, defined by

$$h(t) \stackrel{\text{def}}{=} c(1 + (r-1)\mathbf{1}_{(y,r)}(t))t^{-\alpha}, \quad c > 0, t \in [1, r), \tag{4.2}$$

for any $y \in [1, r)$ leads to the discrete Weibull distribution on $yr^{\mathbb{Z}}$ with parameters $c$ and $\alpha$ and survival function

$$\overline{F}(yr^n) = \exp(-cr^{n\alpha}), \quad n \in \mathbb{Z}.$$

But the ordinary continuous Weibull distribution Weib$(c, \alpha)$ is also covered by the above definition because $h \equiv c$ is an element of $\mathfrak{H}(r, \alpha)$. Let us finally note that

$$\mathcal{L}(W) = r\text{-Weib}(h, \alpha) \iff \mathcal{L}(W^\alpha) = r^\alpha\text{-Exp}(h) \tag{4.3}$$

for all $r, \alpha$ and $h \in \mathfrak{H}(r, \alpha)$, where $h_{1/\alpha}(t) \stackrel{\text{def}}{=} h(t^{1/\alpha}) \in \mathfrak{H}(r^\alpha, 1)$.

We are now ready to give all positive solutions to (1.1) if $T$ possesses a characteristic exponent, i.e. under (A5) or (A6). For convenience we put

$$\text{Exp}(c) = \text{Weib}(c, \alpha) = r\text{-Exp}(h) = r\text{-Weib}(h, \alpha) \stackrel{\text{def}}{=} \delta_0$$

if $c = 0$, respectively $h \equiv 0$.

**Theorem 4.2.** *Let $|J| \geq 2$, $\inf_{j \in J} T_j > 1$ and suppose that $T$ has the characteristic exponent $-\alpha < 0$, i.e. $\sum_{j \in J} T_j^{-\alpha} = 1$.*

(a) *If $\mathbb{G}_\Lambda = \mathbb{R}^{>}$, then*

$$\mathfrak{F}_T^+ = \{\text{Weib}(c, \alpha) \colon c > 0\} \quad [= \{\text{Exp}(c) \colon c > 0\} \text{ if } \alpha = 1].$$

(b) *If $\mathbb{G}_\Lambda = r^{\mathbb{Z}}$ for some $r > 1$, then*

$$\mathfrak{F}_T^+ = \{r\text{-Weib}(h, \alpha) \colon h \in \mathfrak{H}(r, \alpha)\} \quad [= r\text{-Exp}(h) \text{ if } \alpha = 1]$$

*for some $h \in \mathfrak{H}(r, \alpha)$.*

**Proof.** Note first that, for any $\beta > 0$, $T^\beta$ has characteristic exponent $-\beta\alpha$, and $\mathfrak{F}_T = (\mathfrak{F}_{T^\beta})^{1/\beta}$ (see (1.7)). Therefore we can assume without loss of generality that $T$ has characteristic exponent $-1$. Let $F$ be any nontrivial fixed point. It then follows from Lemma 3.2 or Proposition 3.3 that, for some $c > 0$,

$$\nu_F(\mathrm{d}x) = \begin{cases} c\mathbf{1}_{(0,\infty)}(x)\,\mathrm{d}x & \text{if } \mathbb{G}_\Lambda = \mathbb{R}^>, \\ \int_{[1,r)} \sum_{n \in \mathbb{Z}} cr^n \delta_{yr^n}(\mathrm{d}x) \hat{F}(\mathrm{d}y) & \text{if } \mathbb{G}_\Lambda = r^{\mathbb{Z}} \text{ for } r > 1, \end{cases}$$

see (3.5), and thus

$$\nu_F(t) = \begin{cases} ct & \text{if } \mathbb{G}_\Lambda = \mathbb{R}^>, \\ \int_{[1,r)} \sum_{n:\ r^n < t/y} cr^n \hat{F}(\mathrm{d}y) & \text{if } \mathbb{G}_\Lambda = r^{\mathbb{Z}} \text{ for } r > 1, \end{cases}$$

for all $t > 0$. Since $\overline{F}(t) = \exp(-\nu_F(t))$ we infer $F = \mathrm{Exp}(c)$ if $\mathbb{G}_\Lambda = \mathbb{R}^>$. In the discrete case $\mathbb{G}_\Lambda = r^{\mathbb{Z}}$ write $t = r^{m+\beta}$ and $y = r^\gamma$, where $m \in \mathbb{Z}$ and $\beta, \gamma \in [0,1)$. So $m = \lfloor \log_r t \rfloor$ and $\beta = \log_r t - \lfloor \log_r t \rfloor$. Then

$$\nu_F(t) = \int_{[1,r)} \sum_{n < m+\beta-\gamma} cr^n \hat{F}(\mathrm{d}y) = c\left( \int_{[1,r^\beta)} \sum_{n \leq m} r^n \hat{F}(\mathrm{d}y) + \int_{[r^\beta,r)} \sum_{n < m} r^n \hat{F}(\mathrm{d}y) \right)$$

$$= cr^m \hat{F}([1,r^\beta)) + \frac{cr^m}{r-1} = h(t),$$

where

$$h(t) \stackrel{\mathrm{def}}{=} cr^{-\beta}\left( \hat{F}([1,r^\beta)) + \frac{1}{r-1} \right), \quad t > 0,$$

is obviously an element of $\mathfrak{H}(r,1)$. Consequently, $F = r\text{-}\mathrm{Exp}(h)$. Conversely, it is immediately checked that any $r$-periodic exponential distribution with parameter $h \in \mathfrak{H}(r,1)$ is indeed a solution to (1.1). This finishes the proof of the theorem. □

## 5. Solutions for the cases (C2) and (C3)

As announced in the Introduction we will finally treat the problem of finding the solutions of Eq. (1.1) for a vector $T$ containing only negative components $T_j$ (case (C2)), or positive as well as negative ones (case (C3)). In the latter case this is easily accomplished by an appropriate use of the results in the case of positive $T_j$, see Theorem 5.2. However, we begin with a result in the case (C2).

**Theorem 5.1.** *Suppose all $T_j$, $j \in J$, be negative.*

(a) *If $J = \{1, \ldots, n\}$ for some $n \geq 1$, let $\alpha$ be the unique solution of $\alpha + \alpha^n = 1$ in $(0,1)$. Then $\mathfrak{F}_T \setminus \{\delta_0\}$ consists of all $\alpha G + \alpha^n \mathbb{U}_T G$ where $G$ is any distribution on $\mathbb{R}^<$ satisfying*

$$1 - \alpha G(t) = \prod_{i=1}^{n}\left( 1 - \prod_{j=1}^{n} \alpha G\left( \frac{t}{T_i T_j} \right) \right) \tag{5.1}$$

*for $t \leq 0$ and the distribution $\mathbb{U}_T G$ on $\mathbb{R}^>$ is defined by*

$$\mathbb{U}_T G(t) \stackrel{\mathrm{def}}{=} 1 - \prod_{j=1}^{n} G\left( \frac{t}{T_j} + \right), \quad t > 0. \tag{5.2}$$

(b) *If $J = \mathbb{N}$, then $\mathfrak{F}_T = \{\delta_0\}$.*

Given i.i.d. random variables $W, W(1,1), W(1,2), \ldots, W(n,n)$ with common distribution $H \stackrel{\text{def}}{=} \alpha G + (1-\alpha)\delta_0 = \alpha G + \alpha^n \delta_0$, it is easily verified that Eq. (5.1) corresponds to the stochastic fixed point equation

$$W \stackrel{\text{d}}{=} \min_{1 \leq i \leq n} \max_{1 \leq j \leq n} T_i T_j W(i,j). \tag{5.3}$$

As a consequence, an analysis of Eq. (1.1) with negative $T_j$ is very different from that for the case where all $T_j$ are positive. In support of this statement we wish to point out that Theorem 5.1(a), though providing a characterization of the solutions in the case (C2) if $J$ is finite, leaves open the question whether such solutions do exist at all. An answer does indeed require very different techniques from those used here and can be found in the recent article [3]. However, as already mentioned in the Introduction, for the special case $J = \{1, \ldots, n\}$ for some $n \geq 2$ and $T_1 = \cdots = T_n = c < -1$, Khan, Devroye and Neininger [2] could show in the context of randomized game tree evaluation that a nontrivial fixed point exists.

**Proof of Theorem 5.1.** Given negative $T_j$, Eq. (1.1) in terms of the distribution function $F(t) \stackrel{\text{def}}{=} \mathbb{P}(W < t)$ takes the form

$$\overline{F}(t) = \mathbb{P}(T_j W_j \geq t \text{ for all } j \in J)$$
$$= \mathbb{P}\left(W_j \leq \frac{t}{T_j} \text{ for all } j \in J\right) = \prod_{j \in J} F\left(\frac{t}{T_j}+\right), \quad t \in \mathbb{R}. \tag{5.4}$$

If $J$ is finite, then also

$$\overline{F}(t+) = \mathbb{P}\left(W_j < \frac{t}{T_j} \text{ for all } j \in J\right) = \prod_{j \in J} F\left(\frac{t}{T_j}\right), \quad t \in \mathbb{R}. \tag{5.5}$$

Choosing $t = 0$ in (5.4), we see that

$$\overline{F}(0) = F(0+)^{|J|} \tag{5.6}$$

under the usual convention in case $|J| = \infty$ that the right hand side equals 1, if $F(0+) = 1$, and 0, otherwise.

(a) Suppose $J = \{1, \ldots, n\}$ and $\delta_0 \neq F \in \mathfrak{F}_T$. Put $\alpha \stackrel{\text{def}}{=} F(0)$, $\beta \stackrel{\text{def}}{=} \overline{F}(0+)$ and $\gamma \stackrel{\text{def}}{=} F(0+) - F(0) = F(\{0\})$, so $\alpha + \beta + \gamma = 1$. By (5.5), $\beta = \alpha^n$, while $\beta + \gamma = (\alpha + \gamma)^n$ follows from (5.6). Both equations combined give $\gamma = 0$. Consequently, $\alpha + \beta = 1$ and $\beta = 1 - \alpha = \alpha^n$ which shows that $\alpha = \mathbb{P}(W < 0)$ is the unique solution of $\alpha + \alpha^n = 1$ in $(0,1)$.

Next use (5.5) to see that

$$\mathbb{P}(W \geq t | W > 0) = \frac{\overline{F}(t)}{\beta} = \prod_{j=1}^n \frac{F((t/T_j)+)}{\alpha} = \prod_{j=1}^n \mathbb{P}\left(W \leq \frac{t}{T_j} \,\bigg|\, W < 0\right)$$

for all $t \geq 0$. Setting $G \stackrel{\text{def}}{=} \mathbb{P}(W \in \cdot | W < 0)$, we thus have $\mathbb{P}(W \in \cdot | W > 0) = \mathbb{U}_T G$ and therefore

$$F = \alpha \mathbb{P}(W \in \cdot | W < 0) + \beta \mathbb{P}(W \in \cdot | W > 0) = \alpha G + \alpha^n \mathbb{U}_T G$$

as claimed. Next we must prove that $G$ solves Eq. (5.1) for $t \leq 0$. Since $F(t) = \alpha G(t)$ for $t < 0$, we infer from (5.4)

$$1 - \alpha G(t) = \overline{F}(t) = \prod_{i=1}^n F\left(\frac{t}{T_i}+\right) = \prod_{i=1}^n \left(\alpha + \alpha^n \mathbb{U}_T G\left(\frac{t}{T_i}+\right)\right)$$
$$= \prod_{i=1}^n \left(\alpha + \alpha^n\left(1 - \prod_{j=1}^n G\left(\frac{t}{T_i T_j}\right)\right)\right) = \prod_{i=1}^n \left(1 - \prod_{j=1}^n \alpha G\left(\frac{t}{T_i T_j}\right)\right)$$

for all $t < 0$ and thus also for $t = 0$ by left continuity. This shows (5.1).

Conversely, one can easily check that any $F$ of this form with $G$ and $\alpha$ as stated solves Eq. (1.1).

(b) If $J = \mathbb{N}$, it suffices to recall that any $F \in \mathfrak{F}_T$ must satisfy (5.6), now with $|J| = \infty$, which in turn can only hold if $F = \delta_0$ as one can easily check. $\square$

Our final result provides the set of fixed points of (1.1) in the case (C3) where the vector $T$ contains positive as well as negative components $T_j$.

**Theorem 5.2.** *Let $T = (T_j)_{j \in J}$ be a finite or countable family of real numbers such that $J^> \stackrel{\text{def}}{=} \{j\colon T_j > 0\}$ and $J^< \stackrel{\text{def}}{=} \{j\colon T_j < 0\}$ are both nonempty. Then $F \neq \delta_0$ is a solution to Eq. (1.1) i.f.f. $F(0) = 1$ and $F$ solves the very same equation with $T^> \stackrel{\text{def}}{=} (T_j)_{j \in J^>}$ instead of $T$.*

**Proof.** Let $F \in \mathfrak{F}_{T^>}$ and $F(0) = 1$. Let $W$ and $W_j$, $j \in J$, be independent random variables with distribution $F$. Then $W_j < 0$ a.s. implies

$$\inf_{j \in J} T_j W_j = \inf_{j \in J^>} T_j W_j \stackrel{\text{d}}{=} W$$

and thus $F \in \mathfrak{F}_T$.

For the reverse conclusion, we must only show that any $F \in \mathfrak{F}_T \setminus \{\delta_0\}$ satisfies $F(0) = 1$. But (1.1) implies

$$\overline{F}(0) = \mathbb{P}(W \geq 0) \leq \prod_{j \in J^>} \mathbb{P}(W_j \geq 0) \prod_{j \in J^<} \mathbb{P}(W_j \leq 0) = \overline{F}(0)^{|J^>|} F(0+)^{|J^<|}, \tag{5.7}$$

which in view of $|J^<| \geq 1$ can only hold if $\overline{F}(0) = 0$, or $\overline{F}(0) > 0$ and $F(0+) = 1$. To exclude the latter possibility note first that $\overline{F}(0) = 1$ and $F(0+) = 1$ would give $F = \delta_0$. Left with the case $0 < \overline{F}(0) < 1$ and $F(0+) = 1$, we then have $0 < F(\{0\}) = \overline{F}(0) < 1$. Similar to Eq. (2.3) we have here

$$\overline{F}(t) = \prod_{j \in J_n^>} \overline{F}\left(\frac{t}{L(v)}\right) \prod_{j \in J_n^<} F\left(\left(\frac{t}{L(v)}\right)+\right)$$

for all $t \in \mathbb{R}$ and $n \geq 1$, where $J_n^> \stackrel{\text{def}}{=} \{v \in J^n\colon L(v) > 0\}$ and $J_n^<$ is defined accordingly. Use this equation with $n = 2$ to infer in the same way as in (5.7)

$$\overline{F}(0) \leq \overline{F}(0)^{|J_2^>|} F(0+)^{|J_2^<|} \leq \overline{F}(0)^{|J_2^>|}.$$

But $|J_2^>| = |J^> \times J^>| + |J^< \times J^<| \geq 2$ whence the previous inequality yields the contradiction $F(\{0\}) = 0$. $\square$

### Acknowledgment

The authors are very indebted to an anonymous referee for a careful reading of the original manuscript and a number of valuable suggestions that helped improve the presentation of this article.